\documentclass[11pt]{gen-j-l}
\usepackage{amsmath}
\usepackage{amsfonts,amssymb,amsthm}
\usepackage{graphicx}

\newcommand{\thm}[2]{\begin{#1} #2 \end{#1}}

\newcommand{\num}{\mathcal{N\:}}

\begin{document}


\title{Simple curves on surfaces}

\copyrightinfo{1999}{Igor Rivin}

\author{Igor Rivin}

\address{Mathematics department, University of Manchester, Oxford
Road, Manchester}

 \curraddr{Universit\'e Paris-Sud, Orsay and Institut Henri
 Poincar\'e, Paris}

\email{irivin@ma.man.ac.uk}

\dedicatory{To Larry Siebenmann on the occasion of his sixtieth
birthday}

\thanks{I would like to thank Universit\'e Paris-Sud (Orsay) and
the Institut Henri Poincar\'e, and especially its director, Joseph
Oesterl\'e, for hospitality during the writing of this paper}

\date{today}

\begin{abstract}
We study simple closed geodesics on a hyperbolic surface of genus
$g$ with $b$ geodesic boundary components and $c$ cusps. We show
that the number of such geodesics of length at most $L$ is of
order $L^{6g+2b+2c-6}$. This answers a long-standing open
question.
\end{abstract}

 \maketitle

Let $\mathcal{S}$ be a hyperbolic surface of genus $g$ with $c$
cusps and $b$ boundary components. In this paper we study the set
of simple (that is, without self-intersections) closed geodesics
on $\mathcal{S}$. More precisely we study the counting function
$\num(L, \mathcal{S})$ -- the number of simple geodesics of length
no greater than $L$ on the surface $\mathcal{S}$. We show that
there are constants $c_1$ and $c_2$ (depending on only on $S$),
such that
\begin{equation}
\label{asymptotic}
 c_1 L^{6g-6+2b + 2c} \leq \num(L, \mathcal{S}) \leq c_2 L^{6g-6+2b + 2c}.
\end{equation}

The estimate (\ref{asymptotic}) should be put into proper
perspective, and it is with this end that we give the following
historical summary on the study of simple closed curves on
surfaces. This study goes back all the way to the beginning of the
subject of geometry and topology of surfaces (that is, the work of
Henri Poincar\'e and Max Dehn), and some of the subsequent work we
will mention is a developement (if not actually a repetition) of
this work. Being as it may, one line of inquiry has been group
theoretic: suppose $\gamma$ is an element in the fundamental group
of $S$, how do we decide whether or not $\gamma$ is represented by
a simple loop? This question was probably known to Nielsen for the
case of a punctured torus, however, the earliest reference known
to me is the paper of Osborne and Zieschang \cite{osbz}. In the
general case, the first reasonable algorithm for determining
whether an element of $\pi_1(S$) can be represented by a simple
curve was given by H.~Zieschang \cite{z1,z2}, and D.~Chillingworth
\cite{chill1,chill2}, following earlier work of Reinhart
\cite{rein} -- Zieschang's algorithm is primarily group-theoretic,
whilst Reinhart--Chillingworth is more geometric. This work has
been rendered more explicit by Birman and Series \cite{bs0,bs2},
and roughly at the same time Cohen and Lustig \cite{clus,lustig}
had extended the Birman--Series algorithm to determine the mininum
number of intersections between curves representing two homotopy
classes (which includes the self-intersection number of a single
curve as a special case).

Unfortunately, while this algorithmic work is very interesting
(the words in the fundamental group represented by simple closed
curves are a direct generalisation of Sturm sequences  -- these
are precisely the words in the free group on two generators
represented by simple curves on the punctured torus, as remarked
by Birman and Series in their papers), they do not seem to be
usable for estimating the number of such words as a function of
the length of the word (which is, up to constant factors, the same
as the length of the corresponding geodesic).

 This brings us to the the counting
question. It seems that the problem of counting \textit{all}
geodesics of bounded length has been resolved almost entirely, due
to the work of Delsarte, Huber, and Selberg in the constant
curvature case, and Margulis, Bowen, Ruelle, and others in the
variable curvature case. In all cases the estimate is that the
number of geodesics of length bounded above by $L$ is
\textit{asymptotic} to $\exp(h L)/L$, where $h$ is the topological
entropy of the geodesic flow. In particular, $h=1$ for every
finite area hyperbolic surface. Thus, the growth of the number of
all geodesics on all such surface depends (up to first order)
neither on the topology, nor on the actual hyperbolic metric. This
would seem to indicate that the set of all closed geodesics is not
a very geometric object. Now, for simple geodesics, things are a
lot more subtle. For example, for the simplest hyperbolic surface
-- the thrice--punctured sphere -- there are none. After that,
things become more complicated. For the four times punctured
sphere, Beardon, Lehner, and Sheingorn \cite{bls} had shown that
the number of simple geodesics grew at least linearly and at most
quadratically, as a function of length. Since the four-punctured
sphere and the once-punctured torus are essentially the same, this
implies the same estimate for the torus. On the other hand,
Thurston's theories of measured foliations (see \cite{FLP}) and
(dually) projective laminations imply that:

\begin{equation}
\label{masymptotic}
 c_1 L^{6g-6+2b + 2c} \leq \mathcal{M}(L, \mathcal{S}) \leq c_2 L^{6g-6+2b + 2c},
\end{equation}
where $\mathcal{M}(L, S)$  is The number of collections of
pairwise non-intersecting simple closed geodesics of total length
no greater than $L$ on $\mathcal{S}$.

It \textit{is} allowed to take multiple copies of any given curve;
its contribution to the total length is then multiplied by the
multiplicity.

{\bf Notation.} Such a collection of curves will henceforth be
called a multicurve.

Since on a once-punctured torus every multicurve is connected
(though possibly covers itself multiple times), the estimate
(\ref{asymptotic}) for the torus follows, in essence, from the
estimate (\ref{masymptotic}) when $g=1$, and $b+c = 1.$. This also
implies the estimate (\ref{asymptotic}) for the $4$-times punctured
sphere. This estimate only appeared in print in my paper with Greg
McShane \cite{mr1}. In that paper (see also \cite{mr2}) we
actually show a much stronger result: $N(L, S)$ for $S$ a
punctured torus is \textit{asymptotic} to $c_S L^2$,  where the
coefficient $c_S$ depends on the hyperbolic structure, varies
real-analytically over moduli space of tori, and goes to infinity
at infinity of moduli space (hence is not constant).

For $S$ a surface of genus $2$, the estimate (\ref{asymptotic})
follows from the work of Haas and Susskind \cite{haaskind}.

In general, Birman and Series \cite{bs1} have shown that for any
genus, the number of simple curves (actually the number of curves
with a bounded number of self-intersections) grows at most
polynomially, with the exponent depending on the topological type
of the surface. For simple curves their result follows immediately
from the estimate (\ref{masymptotic}), which provides an upper bound
for the number of simple closed curves of bounded length. The
harder part (and the subject of this paper) is proving the lower
bound, It should be noted that this is claimed (indirectly) in the
paper \cite{mrees}. However, the argument there is 
extremely incomplete, and has never been generally accepted.

More recently, Geoff Mess (private communication) has claimed to have
improved the estimate (\ref{masymptotic}) for $M(S, L)$ (the number of
\textit{multicurves} to an
asymptotic result (crudely speaking, showing that one can choose
$c_1$ and $c_2$ arbitrarily close to each other). Furthermore, he
has claimed to be able to show analytic variation of the resulting
constant over moduli space.

Of course, a really interesting question is whether there is an
asymptotic form of (\ref{asymptotic}). The argument proving the
estimate (\ref{asymptotic}) (and occupying the rest of this paper)
seems to indicate that such a result should exist, but it seems
difficult.

Here is an outline of the rest of the paper:

 Section \ref{basics} contains some background facts.
 In Sections  \ref{lowgen} and \ref{genone} the basic method is developed and used to prove estimate
(\ref{asymptotic}) for $g=0, 1$. The case of arbitrary genus
requires a couple of other refinements, and is addressed in
Section \ref{arbgen}.

\noindent {\bf Notation.}In the sequel, whenever constants are
used (denoted by $c$, $c_1$, etc), it is to be understood that
these depend solely on the hyperbolic metric on the surface in
question. The same letter can (and does) denote different numbers
in different places in the paper.

\section{Background}
\label{basics}
 In this section we assemble some necessary
background facts.

\thm{theorem} { \label{cuspneighb} Let $S$ be a hyperbolic
surface, and $\mathcal{C}$ a cusp of $S$. Then there is an
embedded horodisk neighborhood of $\mathcal{C}$ of area $2$ which
contains meets no simple closed geodesic of $S$. }

\begin{proof}
This theorem goes back to Poincar\'e, for a proof see \cite{mr1}.
\end{proof}

\thm{theorem}{\label{twistineq} Let $\gamma$, $\beta$ be simple
closed geodesics on a hyperbolic surface $S$, and let
$\mathcal{T}_\beta(\gamma)$ be the Dehn twist of $\gamma$ around
$\beta$. Then $\ell(\mathcal{T}_\beta(\gamma)) \leq \ell(\gamma) +
i(\gamma, \beta)\ell(\beta)$, where $i(\gamma, \beta)$ is the
geometric intersection number of $\gamma$ and $\beta$. }

\begin{proof} This follows (among other things) from the
variational formula of Wolpert \cite{wol} on the change of length
of curves under earthquake deformations. See also \cite{cp}.
\end{proof}

\thm{theorem} {\label{minintnum} Given two closed curves
$\beta_0$, $\gamma_0$, the smallest geometric intersection number
between curves $\beta$, $\gamma$, freely homotopic to $\beta_0$,
$\gamma_0$, respectively, is realized for $\beta$, $\gamma$
geodesic. }

\begin{proof}
This was also known to Poincar\'e, for a proof see \cite{rlul}.
\end{proof}

From now on, $\Gamma_{A_1\dots A_n}$ shall denote the sphere $S^2$
missing $n$ disks (usually equipped with a hyperbolic metric with
$n$ geodesic boundary components).

\thm{theorem} {\label{class1} There is only one homotopy class of
non-boundary-parallel simple curves on $\Gamma_{ABC}$ beginning and
ending on the same boundary component $A$.}

\begin{proof}
See diagram 1.
\end{proof}

\begin{figure}
\label{onehom}
\includegraphics*[width=2in,height=2in]{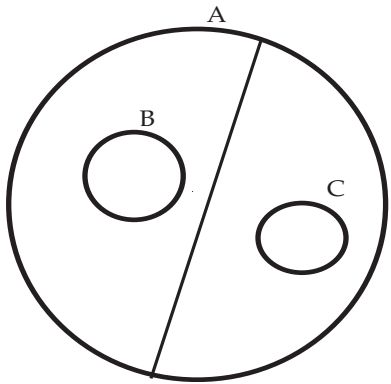}
\caption{first kind}
\end{figure}

\thm{theorem}{\label{class2} There is only one homotopy class of
curves on $\Gamma_{ABC}$ joining boundary component $A$ to
boundary component $B$.}

\begin{proof}
See diagram 2.
\end{proof}

\begin{figure}
\label{twohom}
\includegraphics*[width=2in,height=2in]{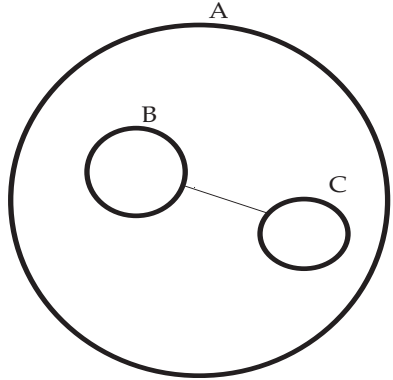}
\caption{second kind}
\end{figure}

Theorem \ref{class1} is relevant in view of:

\thm{theorem} {\label{geodseg} No geodesic segment in
$\Gamma_{ABC}$ can be boundary parallel.}

\begin{proof}
There are no geodesic bigons in $\mathbb{H}^2$.
\end{proof}

These observations are sufficient for us to begin counting simple
curves.

\section{Low genus}
\label{lowgen}

To be systematic, we start with the easiest case:

\thm{theorem} {\label{threeps} Any simple closed geodesic on
$\Gamma_{ABC}$ is a boundary component. }

\begin{proof} Exercise. \end{proof}

\thm{theorem} {\label{fourps}For any hyperbolic structure
$\mathcal{S}$ on the $4$-punctured sphere, $\exists
c_{\mathcal{S}}, L_0 > 0,$, such that for any $L > L_0$, the
number of simple closed geodesics on $\mathcal{S}$ for length not
exceeding $L_0$ is not less than $c_{\mathcal{S}} L^2$.}

\thm{remark}{Of course, a stronger statement follows from
\cite{mr1,mr2}, see the Introduction, but we use the proof of this
result to introduce the techniques and notation for the rest of
the paper.}

\begin{proof}
Let $\Gamma_{ABCD}$ be the $4$-punctured sphere in question. Pick
a simple loop $E$, separating $\Gamma_{ABCD}$ into two
thrice-punctured spheres, $\Gamma_{ABE}$, and $\Gamma_{CDE}$ (see
diagram 3).
\begin{figure}
\label{sepsphere}
\includegraphics*[width=3in,height=2in]{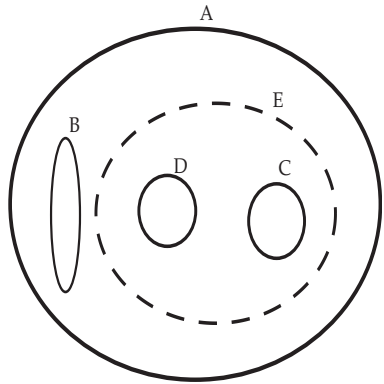}
\caption{A four punctured sphere divided in two}
\end{figure}
Let $\gamma$ be a simple geodesic. There are a finite number of
such which do not intersect the separating curve $E$ (either $4$
or $0$, depending on whether or not one counts boundary curves).
The curve $\gamma$ could be homotopic to $E$, but if not, it must
intersect $E$ transversely in $2k$ points. Let $\gamma$ be such a
curve. Note that $\gamma \cap \Gamma_{ABE}$ consists of $k$
geodesic segments, having all of their endpoints on $E$. Up to
homotopy, there is only one way to thus place $k$ segments in
$\Gamma_{ABE}$, by Theorem \ref{class1}. See Figure 4.
\begin{figure}
\label{ksegs}
\includegraphics*[width=2in,height=2in]{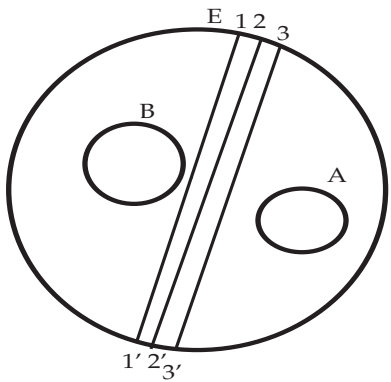}
\caption{$k$ segments from $E$ to $E$}
\end{figure}
The intersection of $\gamma$ with $\Gamma_{CDE}$ looks similar.
Note that the length of $\gamma \cap \Gamma_{CDE} \asymp k.$

Consider the inverse operation: Given two diagrams which look like
figure 4
we can glue them together with a \textit{rational twist}
$\frac{p}{2k}$. The integer part $t=\left\lfloor
\frac{p}{2k}\right\rfloor$ corresponds to twisting $\gamma$ $t$
times around $E$. The fractional part correspond to the change of
the identification map: A diagram (with an orientation) has a
canonical labelling (shown in Figure 4
). Gluing with no
twist corresponds to attaching the strand labelled $1$ to one
labelled $1$, and so on. Twisting by $q$ corresponds to gluing the
strand labelled $l$ to the strand labelled $l+q$, if $l+q < k$, to
$(2k-(l+q))^\prime$ if $l < k$, $l+q > k$ and so on.

It is not hard to see that the twist $\frac{p}{2k}$ leads to a
\textit{connected} curve if and only if $p$ and $2k$ are
relatively prime. By observation \ref{twistineq}, the length of
$\gamma$ twisted by $\frac{p}{2k}$ is bounded above by $p \ell(E)
/ 2$, while the number of twists not exceeding $N$ leading to
connected curves is
\begin{equation}
\frac{\phi(k)}{k} 2 k N = 2 \phi(k) N,
\end{equation}
where $\phi$ denotes the Euler totient function.

By the previous observations, the length of curves obtained
thereby is bounded above by $k N \ell(E)$, so to obtain curves of
length not exceeding $L$, we must take $N \leq \frac{L}{k
\ell(E)}$, thus, for a fixed $k$ we have
\begin{equation*}
\frac{\phi(k) L}{k \ell(E)}
\end{equation*}
curves.
Since $k$
could be anything up to $c L$ (the constant $c$ depending on the
metric on $\Gamma_{ABCD}$), we see that 
\begin{equation}
\label{fourest}
\mathcal{N}(L,\Gamma_{ABCD}) \geq \sum_{k=1}^{c L}
\frac{L}{\ell(E)}\frac{\phi(k)}{k} \geq c^\prime L^2,
\end{equation}
 where we use $\mathcal{N}(L, \mathcal{S})$ to denote the number
 of simple geodesics on a surface $\mathcal{S}$ of length not
 exceeding $L$. The last inequality in Eq. (\ref{fourest}) is the
 consequence of Lemma \ref{euler} below (the argument is standard in
 number theory; we give it here for completeness).
\end{proof}

\thm{lemma}
{
\label{euler}
\begin{equation*}
\sum_{m=1}^n \frac{\phi(m)}{m} = \frac{6}{\pi^2} n + O(\log n),
\end{equation*}
where $\phi$ denotes the Euler totient function.
}

\begin{proof}
First, note that
\begin{equation}
\label{phieq}
\sum_{d | n} \phi(d) = n.
\end{equation}
This follows, for example, from the observation that the number of
elements of order $d | n$ in the cyclic group of order $n$ is equal to
$\phi(d)$. From equation (\ref{phieq}), we have, by M\"obius inversion,
that 
\begin{equation}
\label{mob}
\frac{\phi(n)}{n}=\sum_{d | n} \frac{\mu(d)}{d}.
\end{equation}
Using equation (\ref{mob}) we have
\begin{equation*}
S(n) = \sum_{m\leq n} \frac{\phi(m)}{m} = \sum_{m\leq n} \sum_{d | m}
\frac{\mu(d)}{d}.
\end{equation*}
Changing the order of summation, we see that
\begin{equation*}
S(n) = \sum_{d \leq n}\sum_{j\leq \frac{n}{d}} \frac{\mu(d)}{d} = 
\sum_{d\leq n} \left\lfloor{\frac{n}{d}}\right\rfloor
\frac{\mu(d)}{d},
\end{equation*}
where $\lfloor x \rfloor$ denotes the integer part of $x$. Since $| x
- \lfloor x \rfloor | < 1,$ we have the estimate
$$|S(n) - S_0(n)| < \sum_{d\leq n} |\frac{\mu(d)}{d}| <
\sum_{d<n}\frac{1}{d} = \log(n) + O(1),$$
where 
$$S_0(n) = \sum_{d\leq n} {\frac{n}{d}}
\frac{\mu(d)}{d} = n \sum_{d \leq n} \frac{\mu(d)}{d^2}.$$
Note that 
$$\sum_{j=1}^\infty \frac{\mu(j)}{j^s} = \frac{1}{\zeta(s)},$$
where $\zeta$ is the Riemann $\zeta$ function,
while
$$|\sum_{j=n+1}^\infty \frac{\mu(j)}{j^2}| \leq \sum_{j=n+1}^\infty
\frac{1}{j^2}| = O(n^{-1}).$$ 
Putting all these estimates together, we get the conclusion of the
lemma.
\end{proof}

This analysis will now be extended to deal with a sphere with $k$
boundary components. We will prove

\thm{theorem} {\label{gensph} The number of simple geodesics of
length bounded by $L$ on a sphere with $c$ boundary components
grows like $L^{2c-6}$. }

\begin{proof}
As before, take the sphere $\Gamma_{A_1\dots A_k}$ and cut it into
a sphere $\Gamma_{A_1\dots A_{k-2}E}$ with $k-1$ boundary
components and a pair of pants $\Gamma_{EA_{k-1}A_k}$. We will,
again, count the simple geodesics $\gamma$ which intersect $E$
$2k$ times, then sum over the possible values of $k$. The
intersection of such a curve with $\Gamma_{EA_{k-1}A_k}$ was
already studied in the proof of Theorem \ref{fourps}. It remains
to analyse $\gamma \cap \Gamma_{A_1\dots A_{k-2}E}$. This is a
collection of $k$ disjoint segments having all of their endpoints
on $E$, and we will analyse this inductively. We will prove the
following

\thm{lemma}{\label{onended} Let $G=\Gamma_{A_1\dots
A_{\mathfrak{c}}}$ be a sphere with $\mathfrak{c}$ boundary
components. The number of $k$ component multicurves on $G$ of
length bounded above by $L$ is bounded below by a constant times
$k^{\mathfrak{c}-3}L^{\mathfrak{c}-3}$.}

\begin{proof}[Proof of lemma \ref{onended}] Consider a sphere with $\mathfrak{c}$ boundary components,
$\Gamma_{A_1\dots A_{l-1}E.}$ We cut off a $3$-punctured sphere
$\Gamma_{FA_{l-1}E}$. The intersection of $\gamma$ with this
sphere is a collection of segments, $2 m$ of which go from $F$ to
$E$, while $k-m$ go from $E$ to itself (see diagram 5).
\begin{figure}
\label{segEF}
\includegraphics*[width=2in,height=2in]{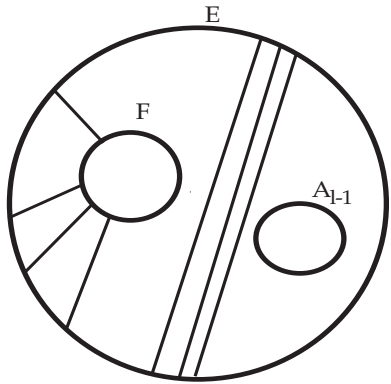}
\caption{$4$ segments from $F$ to $E$, $3$ segments from $E$ to
$E$}
\end{figure}
\thm{remark} {There is another combinatorial possibility, shown in
Figure 6
, but as we are only interested in a lower
bound, we ignore it.} The intersection of $\gamma$ with
$\Gamma_{A_1\dots A_{\mathfrak{c}-1}F}$ has $m<k$ connected
components.
\begin{figure}
\label{segff}
\includegraphics*[width=2in,height=2in]{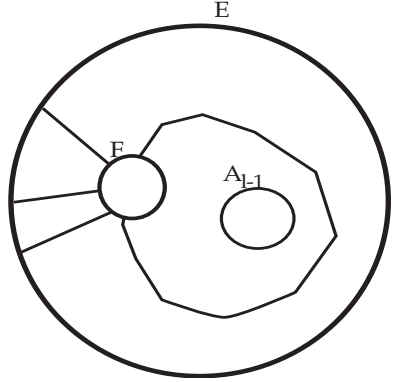}
\caption{Another possibility}
\end{figure}
Thus, we have the inequality
\begin{equation}
\label{fundrec} \mathcal{N}(\mathfrak{c}, k, L) \geq c
\sum_{m=0}^k \int_0^L(L-x) d\mathcal{N}(\mathfrak{c-1}, m, x),
\end{equation}
where the integral is in the sense of Stieltjes, and
$\mathcal{N}(\mathfrak{c}, l, L)$ is the number of $l$-component
multicurves of length bounded above by $L$ on a sphere with
$\mathfrak{c}$ boundary components beginning and ending on a fixed
component.

\thm{remark} {\label{intparts} $f(x)\asymp x^k, k \geq 0$, then
\begin{equation}
\int_0^L (L-x) d f(x) = \int_0^L f(x) d x \asymp x^{k+1},
\end{equation}
integrating by parts.}

\begin{proof}[Proof of inequality (\ref{fundrec})]
Note that each term in the sum comes from the intersection of
$\gamma$ with $\Gamma_{A_1\dots A_{\mathfrak{c}-1}F}$ having $m$
components. In each case, if that intersection has length $x$, we
can twist  the length of that intersection is $x$, and we can
twist a number of times around $F$  to bring the length up to $L$.
This number is proportional to $L-x$, (up to a constant, of order
of length of $F$). \end{proof}

 For example, for $c=3$,
\begin{equation*} \num( 3, l, L) =
  \begin{cases}
   0 & \text{$L < L_0$}, \\
    1 & \text{otherwise}.
  \end{cases}
\end{equation*}
For $c=4$,
\begin{equation}\label{eq4}
\mathcal{N}( 4, l, L) \geq c l L,
\end{equation}
by virtue of Remark \ref{intparts}. An easy inductive argument
shows that
\begin{equation}\label{eegen}
\mathcal{N}(\mathfrak{c}, l, L) \geq c l^{\mathfrak{c}-3}
L^{\mathfrak{c}-3}.
\end{equation}
\end{proof}

To complete the proof of Theorem \ref{gensph}, use Lemma
\ref{onended}, and essentially repeat the counting argument in the
proof of that Lemma \textit{verbatim}.
\end{proof}

\section{Surfaces of genus one}
\label{genone}

To extend our estimates to surfaces of genus $1$, we follow the
same basic strategy: given a surface of genus $1$ with $c$
boundary components (we will denote such a surface by $T_{A_1\dots
A_c}$ , we cut along a nonseparating simple curve $E$ to obtain a
surface of genus $0$ with $c+2$ boundary components: $\Gamma_{E_1
A_1 \dots A_c E_2}.$ Given a simple curve $\gamma$ on $T_{A_1\dots
A_c}$, $\gamma \cap \circ{\Gamma}_{E_1 A_1 \dots A_c E_2}$ is a
collection of segments, each of whose endpoints lies either on
$E_1$ or $E_2$. Conversely, given such a collection of segments,
of total length $x$, we know that by fractional twisting we can
produce several curves $\gamma_1, \dots, \gamma_N$, by identifying
the boundary components $E_1$ and $E_2$ with a fractional twist.
The number $N$ of such curves is of order $L-x$, just as in the
proof of Theorem \ref{fourps}. We now need to count collections of
segments as above, having $k$ intersections with $E_1$ and $k$
intersections with $E_2$. In fact, for the purposes of induction
we will count collections of segments having $k$ intersections
with $E_1$ and $l$ intersections with $E_2$, a quantity we shall
denote by $\num(c, k, l, L)$ ($L$ being the upper bound on the
length). For $c=1$, the integers $k$ and $l$ determine the curve
completely, so
\begin{equation} \num(3, k, L) =
  \begin{cases}
   0 & \text{$L < L_0$}, \\
    1 & \text{otherwise}.
  \end{cases}
\end{equation}

For bigger $c$, we cut the $\Gamma_{E_1 A_1 \dots A_c E_2}$ into
two pieces: $\Gamma_{E_1A_1F}$ and $\Gamma_{FA_2\dots A_c E_2}.$
We now classify the possibilities with respect to the number $m$
of intersections with the separating curve $F$, and to simplify
matters, we will assume that $m \leq l$ --
this will give a sufficiently good lower bound. This restricts the
combinatorial possibilities of the intersection of $\gamma$ with
$\Gamma_{E_1A_1F}$ to one: see figure 7
.
\begin{figure}
\label{mlessl}
\includegraphics*[width=3in,height=2in]{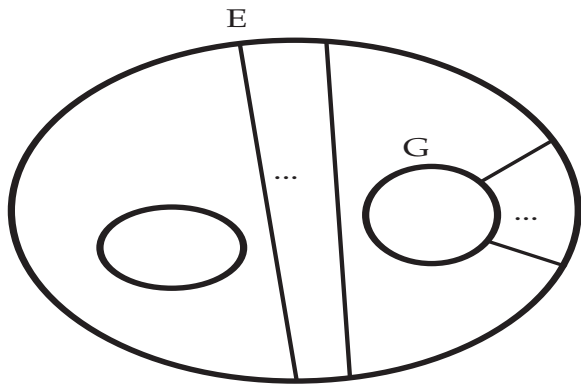}
\caption{Intersection of $\gamma$ with $\Gamma_{E_1A_1F}$}
\end{figure}
 The same argument as in the proof of Theorem
\ref{gensph} gives \begin{equation} \num(\mathfrak{c}, l, k, L)
\geq \sum_{m\leq l} \num(\mathfrak{c}-1, m, k, x) (L-x),
\end{equation}
which gives the following estimate by the same argument as in the
proof of Theorem \ref{gensph}.

\thm{theorem} {\label{tgen1} The number of simple geodesics not
longer than $L$ on a surface of genus $1$ with $c$ punctures is of
order $L^{2c}.$ }

\section{Arbitrary signature}
\label{arbgen}

The main difference between the lower genus situation covered in
the last two sections and the higher genus case considered now is
in the analysis of which fractional twists give connected curves.
In the low genus case, the cyclic order of the intersections of
$k$ segments with a closed loop is always the same: $1, 2, \dots,
k, k^\prime, \dots, 2^\prime, 1^\prime$, which simplifies the
analysis considerably. By contrast, in higher genus, many more
permutations are possible, and it is not \textit{a priori} obvious
how to deal with them. It turns out that we can avoid dealing with
the problem entirely: we only count those curves which behave in a
planar fashion, and these suffice for the lower bound that we
seek. It is at first surprising that we will not have thrown out
the baby with the bathwater, but there is a simple heuristic
explanation: since any collection of simple, pair-wise
non-intersecting, and pair-wise non-isotopic curves contains at
most $3g - 3$ elements, any multicurve falls naturally into (at
most) $3g -3$ subsets, the curves in which are pair-wise parallel,
thus the permutation group action is closer to that of $S_g$ than
of $S_k$ (for a $k$ component multicurve), and thus, if we assume
that every permutation is equally likely, we only lose a constant
factor (in fixed genus $g$). The argument in the rest of this
section is a direct counting argument, which appears to bear out
this heuristic reasoning (which seems difficult to push through
directly).

The actual argument proceeds, as before, by cutting up the surface
into simpler pieces.

Consider a closed surface of genus $g$ (the case of arbitrary
signature will follow by combining the analysis below with the
analysis in Section \ref{lowgen}), and cut it along a curve $E$ to
get two pieces: one, $T_E$, a torus with one boundary component,
the other -- $T_E^{g-1}$ a surface of genus $g-1$ with one
boundary component. The intersection of the simple curve $\gamma$
with $E$ will be a collection of $k$ segments having all of their
endpoints on $E$.

The piece more amenable to analysis is $T_E$. We analyse it by
cutting it further into a thrice-punctured sphere
$\Gamma_{F_1EF_2}$.
There are two combinatorial possibilities for a collection of $k$
mutually non-intersecting segments on this surface: one is shown
in figure 9
, the other in figure 8
.
\begin{figure}
\label{goodsegs}
\includegraphics*[width=3in,height=2in]{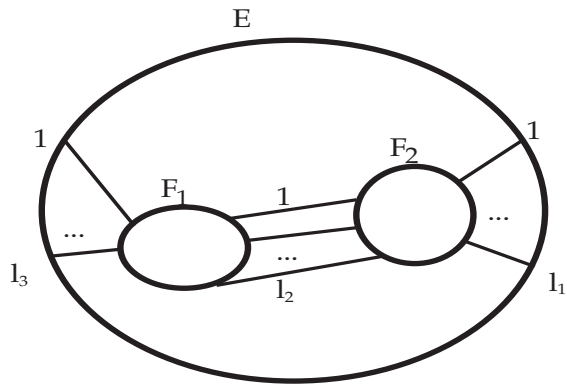}
\caption{Good intersection with a punctured torus}
\end{figure}

\begin{figure}
\label{badsegs}
\includegraphics*[width=3in,height=2in]{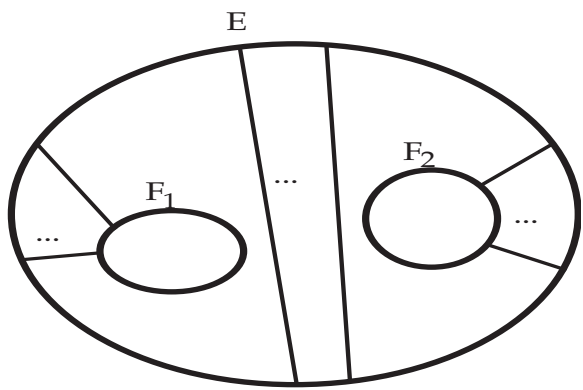}
\caption{Bad intersection with a punctured torus}
\end{figure}

We \textit{forbid} the configuration shown in figure 9
, since this has the wrong cycle type
\textit{vis-\'a-vis} $E$.

We have the constraints that $l_1 + l_3 = k$, while $l_1 + l_2 =
l_3 + l_2$, implying that $l_1 = l_3$.

We are allowed to glue $F_1$ to $F_2$ with a twist, with the
proviso that we get no closed connected components. This seems
like a different sort of problem from the one we encountered in
Sections \ref{lowgen} and \ref{genone}, but luckily there is a
simple trick which allows us to reduce it to that case. To wit, we
glue in a disk with boundary $E$, and use it to connect each
endpoint $i$ to its counterpart $i^\prime$. If, after further
identifying $F_1$ to $F_2$ with a twist, the resulting curve is
connected, obviously no circle components were created. However,
this new problem is exactly the one analysed in Section
\ref{lowgen}. In particular, this tells us that a positive
proportion of the fractional twists are allowed. Now, the total
length of a collection as a function of $l_1$, $l_3$, and the
number $\tau$ of twists is bounded above by $c_1 k + c_2 l_3 + c_3
\tau$, and thus the total number of systems of length bounded
above by $L$ is at least $c (L - c_1 k).$

To analyze the other piece, $T_E^g$, we use an inductive
decomposition, much as before. The case of $g=1$ was done above.
If $g > 1$, we cut the surface into a torus $T_{EF}$ with two
boundary components, and a $T_F^{g-1}$. We are now reduced to
analysing the torus $T_{EF}$. We want information about the
collections of segments which intersect $E$ $k$ times, intersect
$F$ $l$ times, and have a total of $\frac{k+l}{2}$ connected
components. In order to do this we (yet again) cut $T_{EF}$ into
two thrice-punctured spheres $\Gamma_{EC_1C_2}$ and
$\Gamma_{FC_1C_2}$. Let the number of intersections with $C_1$ be
$m$, while the number of intersections with $C_2$ be $n$. By our
requirement on the permutation type of intersection with $E$, the
intersection of our system with $\Gamma_{EC_1C_2}$ must look like
Figure 10
, and similarly for $\Gamma_{FC_1C_2}$ (Figure 11
. We know, furthermore, that $k_2 + k_3 = k$, $k_1 + k_2 = m$,
$k_1 + k_3 = n$, which implies that $$k_1 = \frac{m+n-k}{2},\qquad
k_2 = \frac{k+m-n}{2},\qquad k_3=\frac{k+n-m}{2}.$$ We now have to
worry about two problems. One is that certain ``horizontal''
strands (those between $C_1$ and $C_2$) might close up into loops.
The second is that certain strands might enter and leave through
the boundary component $F$. Either way, we would get a
non-connected multicurve. We deal with both of these by, first,
imposing two additional inequalities. $2 k_1 \leq l_1$, and $2 l_1
\leq \min(m,n)$. We then allow only those twists which connect one
of the $k_1$ strands to one of the $l_1$ strands (there are at
least $l_1 -k_1$ such), and in addition imposing the condition
that the horizontal loops do not close up -- this is not extremely
restrictive, since there we can twist by at least $k_1/2$ around
$C_1$, and also by $k_1/2$ around $C_2$. (see figures 10
 and
  for notation).

\begin{figure}
\label{tophalf}
\includegraphics*[width=3in,height=2in]{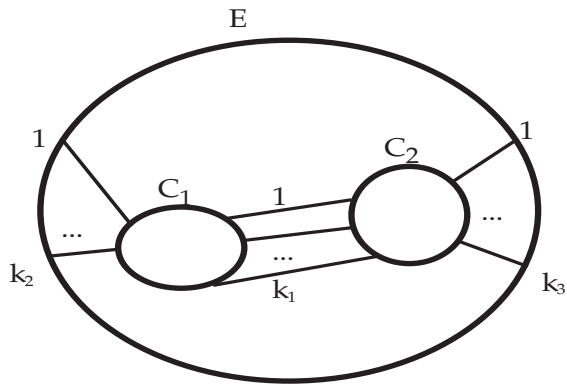}
\caption{Top half of the twice-punctured torus}
\end{figure}

\begin{figure}
\label{bottomhalf}
\includegraphics*[width=3in,height=2in]{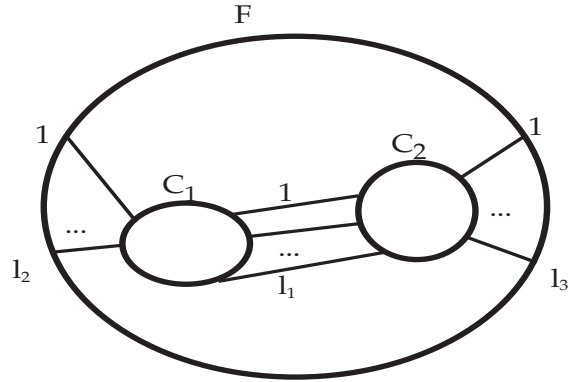}
\caption{Bottom half of the twice-punctured torus}
\end{figure}

It is not hard to see that, subject to these restrictions, any
strand coming into $C_1$ from below will either leave straightaway
upwards, or will cycle around between $C_1$ and $C_2$ for a while
before leaving in that same manner. Likewise for strands coming
into $C_2$ from below. The same sort of counting as before will
show that the set of permissible multicurves of length bounded
above by $x$ is at least $c x^5$, since these corresponds to
points in a cone in $5$ dimensions (corresponding to twisting
around $c_1$, twisting around $c_2$, and the parameters $l, m,
n$). Our constraints are all inequality constraints, and thus will
cut out a non-degenerate cone. The rest of the inductive argument
is as in Sections \ref{lowgen} and \ref{genone}, and so for
compact surfaces we obtain the claimed result

\thm{theorem} {\label{comthm} The number of simple geodesics of
length bounded by $L$ on a compact surface of genus $g$ grows like
$L^{6g-6}$}

The estimate for arbitrary signature follows as indicated in the
beginning of this section.

\section{Conclusions and musings}
The reader will have noted that the estimates on the
\textit{density} of connected multicurves among all multicurves
become worse and worse as the genus of our surfaces increases.
This is, to an extent, a reflection of reality (in fact, it is
easy to see that this density decreases exponentially as a
function of $g+c$. However, it is clear that the estimates one
might obtain by our methods are far from optimal.

Another observation one might make is that the methods of this
paper are obviously insufficient for deriving asymptotic results
(extending those for the punctured torus as in \cite{mr1,mr2}).
While there is some possibility that the method used in the low
genus case might be pushed to get results of this type, for
general genus this seems hopeless. However, the very existence of
an asymptotic formula is in some doubt. A problem which might be
tractable by the current methods is one of finding order of growth
results for curves with a bounded number of self-intersections
(the estimates of Birman and Series are acknowledged by the
authors not to be sharp).

 \specialsection*{Acknowledgements} I would like to
thank Nalini Anantharaman for enlightening discussions,
Caroline Series for help with the historical background, and Greg
McShane for comments on a previous version of this paper.



\bibliographystyle{amsplain}

\begin{thebibliography}{99}

\bibitem{bls}
A.~Beardon, J.~Lehner, M.~Sheingorn.\textit{Closed geodesics on a
Riemann surface with applications to the Markov spectrum,\/} {\em
Transactions of the American Mathematical Society}, {\bf
295}(1986), no. 2, pp.~635--647.


\bibitem{bs0}
Joan Birman and Caroline Series. \textit{An algorithm for simple
curves on surfaces,\/} {\em Journal of the London Mathematical
Society}(2), {\bf 29}(1984), no. 2, 331--342.

\bibitem{bs1}
Joan Birman and Caroline Series. \textit{Geodesics with bounded
intersection number on surfaces are sparsely distributed,\/} {\em
Topology}, {\bf 24}(1985), no. 2, pp.~217-225.

\bibitem{bs2}
Joan Birman and Caroline Series.\textit{Dehn's algorithm
revisited, with applications to simple curves on surfaces,\/} {\em
Combinatorial group theory and topology} (Alta, Utah, 1984),
pp.~451-478, Annals of Mathematics Studies 111, Princeton
University Press, Princeton, NJ, 1987.

\bibitem{chill1}
David Chillingworth. \textit{Winding numbers on surfaces. I,\/}
{\em Math. Ann.} {\bf 196}(1972), pp.~218--249.

\bibitem{chill2}
David Chillingworth. \textit{Winding numbers on surfaces. II,\/}
{\em Math. Ann.} {\bf 199}(1972), pp.~131--153.

\bibitem{clus}
Marshall Cohen and Martin Lustig. \textit{Paths of geodesics and
geometric intersection numbers. I,\/} {\em Combinatorial group
theory and topology} (Alta, Utah, 1984), pp.~479-500, Annals of
Mathematics Studies 111, Princeton University Press, Princeton,
NJ, 1987.

\bibitem{FLP}
Albert Fathi, François Laudenbach, Valentin Poenaru.
\textit{Travaux de Thurston sur les surfaces,\/} With an English
{\em Astérisque}, {\bf 66-67}. Société Mathématique de France,
Paris, 1979..

\bibitem{haaskind}
Andrew Haas and Perry Susskind. \textit{The connectivity of
multicurves determined by integral weight train tracks,\/} {\em
Transactions of the American Mathematical Society,} {\bf
329}(1992), no.~2, pp.~637--652.

\bibitem{lustig}
Martin Lustig. \textit{Paths of geodesics and geometric
intersection numbers. II,\/}{\em Combinatorial group theory and
topology} (Alta, Utah, 1984), pp.~501--543, Annals of Mathematics
Studies 111, Princeton University Press, Princeton, NJ, 1987.


\bibitem{mr1}
Greg McShane and Igor Rivin. \textit{Simple curves on hyperbolic
tori,\/} {\em C.~R.~Acad.~Sci.~Paris S\'er. I. Math.}, {\bf 320},
no. 12, June 1995.

\bibitem{mr2}
Greg McShane and Igor Rivin. \textit{Geometry of geodesics and a
norm on homology,\/} {\em International Mathematics Research
Notices}, February 1995.

\bibitem{messy}
Geoffrey Mess. Private communication.

\bibitem{osbz}
R.~Osborne and H.~Zieschang. \textit{Primitives in the free group
on two generators,\/} {\em Inventiones Mathematicae} {\bf
63}(1981), no. 1, pp.~17--24.


\bibitem{mrees}
Mary Rees. \textit{An alternative approach to the ergodic theory
of measured foliations on surfaces,\/} {\em Ergodic Theory
Dynamical Systems}, {\bf 1} (1981), no. 4, pp.~461--488.


\bibitem{rein}
Bruce Reinhart. \textit{The winding number on two manifolds,\/}
{\em Ann. Inst. Fourier. Grenoble}, {\bf 10}(1960), pp.~271--283.

\bibitem{rlul}
Igor Rivin. \textit{Intrinsic geometry of convex ideal polyhedra
in hyperbolic $3$-space,\/} {\em Analysis, algebra, and computers
in mathematical research (Lule\aa, 1992)}, Lecture Notes in Pure
and Appl. Math., 156, Dekker, New York, 1994, pp.~275-291.

\bibitem{cp}
Caroline Series. \textit{Wolpert's formula,\/} Warwick University
Preprint, 1996.

\bibitem{wol}
Scott Wolpert. \textit{On the symplectic geometry of deformations
of a hyperbolic surface,\/} {\em Annals of Mathematics (2)} {\bf
117}(1983), no.~2, pp.~207--234.

\bibitem{z1}
Heiner Zieschang. \textit{Algorithmen f\"ur einfache Kurven auf
Fl\"achen,\/} {\em Math. Scand.} {\bf 17}(1965), pp.~17--40.

\bibitem{z2}
Heiner Zieschang. \textit{Algorithmen f\"ur einfache Kurven auf
Fl\"achen. II,\/} {\em Math. Scand.} {\bf 25}(1969), pp.~49--58.



\end{thebibliography}

\end{document}